
\input amstex
\documentstyle{amsppt}
\define\longdashrightarrow{{\,\text{- - - -}{\scriptscriptstyle{}^{{}_{>}}}}}
\define\longdashleftarrow{{{\scriptscriptstyle{}^{{}_{<}}}\text{- - - -}\,}}
\define\CDdashright#1#2{&\,\mathop{\longdashrightarrow}\limits^{#1}_{#2}\,&}
\define\CDdashleft#1#2{&\,\mathop{\longdashleftarrow}\limits^{#1}_{#2}\,&}
\define\CDsubset{&\quad\mathop{\subset}\quad&}

\def\P{\Bbb P}

\def\Til#1{\widetilde{#1}}

\def\PGL{\text{PGL}}
\hcorrection{.7 in}

\def\phy{\varphi}

\pagewidth{12.45 cm}\pageheight{20.35 cm}
\topmatter
\title Linear orbits of $d$-tuples of points in $\P^1$
\endtitle
\author Paolo Aluffi \\ Carel Faber\endauthor
\address {Mathematics Department, Florida State University, Tallahassee FL
32306, USA}
\endaddress
\address {Faculteit Wiskunde en Informatica, Univ. van Amsterdam, Plantage
Muidergracht 24, 1018 TV Amsterdam, NL}
\endaddress
\endtopmatter
\document

\heading \S 0. Introduction \endheading

The group $\PGL(2)$ of linear transformations of the projective line
$\P^1$ acts naturally on the set of configurations of points on the
line. We call each configuration of $d$ points (some of which may
coincide in the same point on the line) a \lq $d$-tuple' of points;
for a given $d$, the set of $d$-tuples of points in $\P^1$ forms a
dimension-$d$ projective space $\P^d$. In this note we are concerned
with the orbits of this action of $\PGL(2)$ on $\P^d$. The closure of
each orbit is a projective subvariety of $\P^d$ of which we determine
the degree (\S 1), the \lq boundary'--i.e., the complement of an orbit
in its closure--(\S 2), and the multiplicity at points of the boundary
(\S 3). These results are used to provide a complete classification of
the non-singular orbit closures, and criteria for an orbit closure to
be non-singular in codimension 1 (\S 4).

Although seemingly natural objects of study, we didn't find a lot of
work on these orbits in the literature. Some of the results presented
here appear also in \cite{Mukai-Umemura}, in one form or another; and
the `combinatorial' computation of the degree that we will sketch in
this introduction goes back to \cite{Enriques-Fano}. But for example
Mukai and Umemura establish the non-singularity of the orbit closures
of a specific 6-tuple and a specific 12-tuple by a rather long and
ad-hoc coordinate computation.  We hope to provide here a more
unifying approach.

Our main motivation in this study is to prepare the ground for the
much richer case of the action of $\PGL(3)$ on spaces parametrizing
plane curves. The approach we use in this note is susceptible to be
employed in higher dimensions, although the technical difficulties
mount very rapidly. The reader wishing to approach the $\PGL(3)$ case
(see \cite{Aluffi-Faber}) will find here a sample of the essential
techniques.\vskip 12pt

The main idea for the degree and multiplicity computations is the
following: for each given $d$-tuple of points on $\P^1$, build a
smooth variety $\Til V$ and a proper map from this to the
closure of the orbit of the $d$-tuple. In fact this $\Til V$ will be a
compactification of $\PGL(2)$, determined by the $d$-tuple, which we
obtain by a suitable blow-up of the $\P^3$ of $2\times 2$
(homogeneous) matrices. After the construction, we reduce the
calculations to calculations on $\Til V$, where some intersection
calculus (particularly, the formalism of Segre classes of
\cite{Fulton}) allows us to perform them. The blow-up construction
also allows us to determine explicitly the boundary of the orbit.

The classification of smooth orbit closures follows from the
multiplicity computations of \S 3; we use the classification of finite
subgroups of $\PGL(2)$, which can be found for example in \cite{Weber}.

We now sketch here the easy \lq combinatorial' computation of the
degree of the orbit closure of a $d$-tuple consisting of $d\ge 3$ {\sl
distinct\/} points.  In this case the orbit closure is 3-dimensional,
so its degree may be computed as the intersection product with three
hyperplanes of $\P^d$.

For the hyperplanes, take 3 distinct `point-conditions', i.e.,
hyperplanes in $\P^d$ consisting of the $d$@-tuples that contain a
certain given point. One checks easily that the intersection
multiplicity of the orbit closure and three point-conditions
(determined by three distinct points $p_1$, $p_2$, $p_3$) at a
$d$-tuple equals the product of the multiplicities of $p_1$, $p_2$ and
$p_3$ in the $d$-tuple: so the intersection is automatically {\sl
transversal\/} if the $d$-tuple consists of $d$ {\sl distinct\/}
points. Therefore, in this case the degree is just the number of
points of intersection: the computation then comes down to counting
the number of elements of $\PGL(2)$ that send a given $d$@-tuple
(consisting of $d$ distinct points) to a $d$@-tuple that contains 3
(distinct) given points. Since an element of $\PGL(2)$ is uniquely
determined by prescribing the images of 3 distinct points, one sees
that the answer must be
$$d(d-1)(d-2).$$
To get the degree of the orbit closure, we have to divide this number
by the number of elements of PGL(2) sending a $d$-tuple to itself:
i.e., the order of the stabilizer of the $d$@-tuple. For example:
\roster
\item The stabilizer of a 3@-tuple consisting of 3 distinct points is
$S_3$, so the degree of the orbit closure is 1 (the orbit closure is
$\P^3$).
\item A general 4@-tuple has stabilizer $C_2 \times C_2$, so the
degree of the orbit closure is $\frac {4\cdot 3\cdot 2}4=6$. The
4@-tuples with $j=0$ (resp.~1728) have stabilizers $A_4$
(resp.~$D_4$), so that the orbit closure has degree 2 (resp.~3).
\item For $d\ge5$, a general $d$@-tuple has trivial stabilizer, so the
degree of the orbit closure is $d(d-1)(d-2)$.
\endroster

It would be easy to apply the same procedure to examine the case in
which some points of the $d$-tuples appear with multiplicity. However,
we don't see how to obtain by this approach a unified treatment of all
cases; more importantly, this approach wouldn't help us to study the
singularity of these orbit closures, and more important still we don't
see how this kind of computations could be interpreted to attack
higher dimensional cases such as the one dealt with in
\cite{Aluffi-Faber}.\vskip 12pt

{\sl Acknowledgement.\/} Both authors wish to thank the
Max-Planck-Institut f\"ur Mathematik for the wonderful hospitality.

\heading \S 1. The predegree of the orbit closure. \endheading

We work over an algebraically closed field of characteristic 0.

The first question we consider is the computation of the degree of the
closure (in~$\P^d$) of the orbit of a $d$@-tuple under the action of
$\PGL(2)$. Here we think of $\P^d$ as the space parametrizing
homogeneous forms of degree $d$ on $\P^1$, and each point of this space
is identified with the $d$-tuple of zeros of the form corresponding to
it.  Also, we will denote by $s$ the number of {\sl distinct\/} points
in the $d$-tuple. As mentioned in the introduction, the main
ingredient in the computation is the construction for each $d$-tuple
of a non-singular variety dominating the orbit closure.

First we observe this is not necessary if the whole $d$-tuple is
concentrated in one point (that is, if $s=1$). We'll refer to this
particular $d$-tuple as to the \lq $d$-fold point', and the reader
should have no difficulties in checking that the orbit of the $d$-fold
point (that is, the set of all such $d$-tuples) is simply the
degree-$d$ rational normal curve in $\P^d$.

Next, let's consider the case when the $d$-tuple is distributed among
2 distinct points, that is one $r$-fold point and one distinct
$(d-r)$-fold point. Again, in this case the reader will see
immediately that the orbit consists of all $d$-tuples with the same
multiplicity data.

\proclaim{Proposition 1.1} The orbit closures of $d$-tuples consisting of an
$r$-fold point and a $(d-r)$-fold (distinct) point are surfaces in $\P^d$,
of degree: $2r(d-r)$ if $r\ne d/2$, $r(d-r)=r^2$ if $r=d/2$.\endproclaim

\demo{Proof} For this, we dominate the orbit closure with
$\P^1\times\P^1$, using the map $\P^1\times\P^1 @>>> \P^d$ defined by
$$((a_0:a_1),(b_0:b_1)) \mapsto (a_1 x - a_0 y)^r(b_1 x - b_0 y)^{d-r}
\quad:$$
it is clear that this map is finite, and that the complement of the
diagonal in $\P^1\times \P^1$ maps onto the orbit we are considering.
Also, it is clear that the degree of this map is 1 if $d\ne 2r$, and 2 if
$d = 2r$: so to get the statement we just need to check that the
self-intersection of the pull-back of the hyperplane class from $\P^d$ to
$\P^1\times\P^1$ via the above map is $2r(d-r)$. This is straightforward:
if $h_1$, $h_2$ denote the hyperplane class of the factors, the pull-back
of the hyperplane class from $\P^d$ is $(r h_1+(d-r) h_2)$, and
$$\int_{\P^1\times \P^1} (r h_1+(d-r) h_2)^2 = \int_{\P^1\times \P^1}
2r(d-r) h_1h_2= 2r(d-r)\quad.$$
(Here and in the following $\int$ will denote \lq degree' in the sense of
\cite{Fulton})\qed\enddemo

It's worth observing that if $r=d/2$, then the orbit closure is a
(regular) projection to $\P^d$ of the $r$-th Veronese embedding of
$\P^2$---the degree is indeed $r^2$ in this case, as it should be. For
example, for $r=2$ this is the (non-singular) projection of the
Veronese surface in $\P^5$ to $\P^4$.\vskip 12pt

Now we move to the most interesting case, that of a $d$-tuple distributed
in $s\ge 3$ points. In this case the orbit and its closure have dimension
3. In order to construct a non-singular threefold dominating the orbit
closure of a given $d$-tuple, we resolve the indeterminacies of a rational
map associated naturally to the given $d$-tuple.

Choose coordinates $(x:y)$ in $\P^1$, and let $C$ stand for a
homogeneous form in $(x:y)$ of degree $d\ge3$, and for the $d$@-tuple
of points on $\P^1$ corresponding to it. The $\PGL(2)$@-orbit of $C$
in $\P^d$ is the image of the map
$$c\: \PGL(2) \to \P^d$$
sending $\alpha \in\PGL(2)$ to the form $C \circ \alpha$. Observe that
this map is finite (if at least three points of the $d$-tuple are
distinct), and its degree equals the order of the stabilizer of $C$.
This map determines a rational map from the $\P^3$ of $2\times 2$
matrices to $\P^d$, which we also denote by $c$.

Now we will resolve this rational map: i.e., we will construct a
variety $\Til V$ filling a commutative diagram
$$\CD
\PGL(2) \CDsubset \Til V @>\Til c>> \P^d \\
@| @V\pi VV @| \\
\PGL(2) \CDsubset \P^3 \CDdashright c{} \P^d
\endCD$$

The image of $\tilde c$ in $\P^d$ is precisely the orbit closure. Thus
the degree of the orbit closure can be found by computing the third
power of the pull-back of the hyperplane class of $\P^d$ to $\Til V$,
and dividing by the order of the stabilizer of $C$. We call
`predegree' the product of the degree by the order of the stabilizer:
since the $d$-tuple is supported on at least 3 points, this term will
be synonymous for the 3-fold self-intersection of the pull-back of the
hyperplane from $\P^d$.

The base locus of $c\: \P^3 \dashrightarrow \P^d$ consists of the
matrices $\alpha$ for which the form $C \circ \alpha$ is identically
zero. This happens exactly when $\alpha$ is a rank@-1 matrix with image
a point of the $d$@-tuple $C$. The base locus of $c$ is therefore
supported on a finite number of `parallel' lines in the (non-singular)
quadric of rank@-1 matrices. There are as many distinct lines as there
are distinct points in the $d$@-tuple $C$.

\proclaim{Proposition 1.2} A variety $\Til V$ as above can be obtained
by blowing up $\P^3$ along the support of the base locus of
$c$.\endproclaim

\demo{Proof} To see this, call `point-conditions in $\P^3$' the
inverse image of the point-conditions of $\P^d$ (defined above). The
map $c$ is then the map defined by the linear system generated by the
point-conditions in $\P^3$, and therefore the base locus of $c$ is
actually cut out by the point-conditions. Now we argue that a
point-condition in $\P^3$ is a degree-$d$ hypersurface consisting of
nothing but a collection of hyperplanes, one for each point in the
$d$-tuple $C$, each appearing with the same multiplicity as the
corresponding point appears in~$C$.  This is immediate: give
coordinates

$$
\left( \matrix
p_0 & p_1 \\
p_2 & p_3
\endmatrix \right)
$$
to the $\P^3$ of matrices; and suppose $C$ is given by the equation
$$F(x:y)=0\quad.$$
Then the point-condition corresponding to e.g\. the point $(1:0)$ has
equation
$$F(p_0:p_2)=0\quad,$$
so is indeed a union of hyperplanes as argued.

Let $\Til V$ be the blow-up of $\P^3$ along the lines supporting the
base locus of $c$. The (a~priori rational) map $\tilde c$ making the
above diagram commute is then defined by the linear system on $\Til V$
generated by the proper transforms of the point-conditions: so the
base locus of $\tilde c$ is cut out by the proper transforms in $\Til
V$ of the point-conditions. But since the point-conditions are
supported on unions of {\sl hyperplanes\/}, they necessarily intersect
transversally in $\P^3$ along the base locus of $c$: therefore their
intersection in $\Til V$ is empty, and we can conclude that the map
$\tilde c: \Til V @>>> \P^d$ is indeed a morphism.\qed\enddemo
\vskip 12pt

Now computing the 3-fold self-intersection of the class of the proper
transform of a point-condition (i.e., the predegree of the orbit
closure) is a straightforward intersection calculus exercise. We use
\cite{Aluffi-Faber}, Proposition 3.2: the self-intersection is
computed as the self-intersection of the point-condition in $\P^3$
(i.e., $d^3$) {\sl minus\/} contributions coming from each component
of the base locus of $c$. The formula gives
$$\text{predegree}=d^3-\sum_{i=1}^s \int_{L_i}
\frac{(m_i+dh)^3}{1+2h}\quad,$$
where the summation runs over the {\sl distinct\/} points $p_1,\dots,
p_s$ of the $d$-tuple, $L_i$ is the line in the base locus
corresponding to $p_i$, $m_i$ is the multiplicity of $p_i$ in the
$d$-tuple (thus the multiplicity of the point-conditions along $L_i$),
and $h$ denotes the hyperplane class in $L_i$. The degree is computed
by taking the coefficient of $h$ in the expression under $\int$. Doing
this gives:

\proclaim{Proposition 1.3} For $d\ge 3$, the predegree of the orbit
closure of a $d$-tuple is
$$
d^3-3d(\sum_{i=1}^s m_i^2)+2(\sum_{i=1}^s m_i^3)\quad.
$$
\endproclaim

So the predegree of a $d$-tuple $C$ can be written in terms of just
$d$ and {\sl two\/} numbers, each of which is a sum of `local
contributions' given by each point of $C$. For example, if the
$d$-tuple consists of $d-r$ simple points and one $r$-fold point,
then
$$\sum_{i=1}^s m_i^2=r^2+d-r,\qquad \sum_{i=1}^s m_i^3=r^3+d-r,$$
so
$$\align
\text{predegree}&=d^3-3d(r^2+d-r)+2(r^3+d-r)\\
&=(d-r)(d-r-1)(d+2r-2)\quad.
\endalign$$
As seen in \cite{Aluffi-Faber}, this general feature of the predegree
(being determined by a few numbers recording local data) is preserved
in the $\PGL(3)$ case, at least for smooth curves.

For $s=1$ or 2, the formula of this proposition gives 0: which
reflects the fact that in these cases the orbits have dimension $<3$.
We also remark that the $\P^1\times\P^1$ used to dominate the orbit
closure in the case $s=2$ in Proposition (1.1)  can
also be seen as one component of the exceptional divisor of the same
blow-up construction used for the case $s\ge 3$.\vskip 12pt

\heading \S 2. The boundary of an orbit closure\endheading

We turn now to the question of determining the `boundary' of the orbit
of a $d$-tuple $C$, by which we mean the complement of the orbit in
its closure. Observe that the boundary of an orbit is necessarily
itself the union of orbits, and has dimension $\le 2$. Since the orbit
of a $d$-tuple has dimension 3 as soon as the $d$-tuple consists of at
least 3 distinct points, we can conclude right away that the boundary
of the orbit of a given $d$-tuple must consist of a union of
orbits of $d$-tuples concentrated in at most {\sl two\/} points. We
will show:

\proclaim{Proposition 2.1} The boundary of the (3@-dimensional) orbit
of $C$ is the union of the 1@-dimensional orbit of $x^d$ and of those
2@-dimensional orbits of $x^ry^{d-r}$ for which $r$ is the
multiplicity of a point of $C$.\endproclaim

\demo{Proof} We use again the variety $\Til V$ constructed in \S1. The
rank-1 matrices {\sl not\/} in the base locus have image in the orbit
of $x^d$; so we only have to determine the image in $\P^d$ of the
components of the exceptional divisor in $\Til V$. Give coordinates
$$
\left( \matrix
p_0 & p_1 \\
p_2 & p_3
\endmatrix \right)
$$
to the $\P^3$ of matrices; the locus of rank@-1 matrices is given by
$p_0p_3-p_1p_2=0$. Suppose the $d$-tuple $C$ has equation
$a_0x^d+a_1x^{d-1}y+\dots+a_dy^d=0$, corresponding to the point
$(a_0:a_1:\dots:a_d)\in \P^d$ (with obvious choice of coordinates
there). Assume that $(1:0)$ is a point of multiplicity $r\ge1$ in $C$,
i.e., $a_0=a_1=\dots=a_{r-1}=0, a_r\ne 0$. Then $p_2=p_3=0$ is a
component of the base locus of $c$ and we can study $\Til V$ locally
by blowing up $\P^3$ along $p_2=p_3=0$.

On the affine piece $p_0=1$ we have coordinates $(p_1, p_2, p_3)$. On an
affine piece of the blow-up, coordinates $(q_1, q_2, q_3)$ are given by
$$\left\{
\aligned
p_1&=q_1 \\
p_2&=q_2 \\
p_3&=q_2q_3
\endaligned \right.
$$
The map induced by $c$ is then given by
$$ (q_1, q_2, q_3) \mapsto (b_0:b_1:\dots:b_d)$$
with
$$
b_0x^d+\dots+b_dy^d \sim
a_r(x+q_1y)^{d-r}(q_2x+q_2q_3y)^r+\dots+a_d(q_2x+q_2q_3y)^d.
$$
Note that we can factor out $q_2{}^r$ from the last expression, so that
$$\multline
b_0x^d+\dots+b_dy^d \sim a_r(x+q_1y)^{d-r}(x+q_3y)^r \\
+a_{r+1}q_2(x+q_1y)^{d-r-1}(x+q_3y)^{r+1}+\dots+a_dq_2^{d-r}(x+q_3y)^d.
\endmultline
$$
The exceptional divisor is given here by $q_2=0$. The restriction of
the map $\tilde c: \Til V @>>> \P^d$ to the component of the
exceptional divisor of $\Til V$ corresponding to the $r$-fold point is
then given by restricting the last expression to $q_2=0$: we get
$d$-tuples corresponding to points
$$b_0x^d+\dots+b_dy^d \sim a_r(x+q_1y)^{d-r}(x+q_3y)^r:\tag{*}$$
we conclude that the image of the exceptional divisor corresponding to
a point in $C$ of multiplicity $r$ is the closure of the
$\PGL(2)$@-orbit of $x^{d-r}y^r$. (The boundary of {\sl this\/} orbit
is the orbit of $x^d$.) The statement follows.\qed\enddemo\vskip 12pt

\heading \S 3. Multiplicities. \endheading

We will now use the blow-up construction described in \S 1 to compute
the multiplicity of the closure of an orbit along the orbits making up
its boundary. For $s=1$ and $s=2,\, r=d/2$ (notations as in \S 1) we
have remarked that the orbit closure is essentially a Veronese, so it
is non-singular. To analyze the situation for $s=2,\, r\ne d/2$ and
$s\ge 3$, we first need the following fact.

Identify $\P^d$ with the space of $d$-tuples of points on $\P^1$, by
giving it coordinates $(a_0:\dots:a_d)$ and associating with every
$A=(a_0:\dots:a_d)$ the $d$-tuple of zeros of
$F_A(x:y)=a_0x^d+a_1x^{d-1}y+\dots+a_dy^d$. Then let $H_A(x:y)$ denote
the Hessian of this form with respect to $x,y$, a form itself of
degree $2d-4$ in $(x:y)$ for each given $A$. For a given $(\xi:\eta)$
in $\P^1$, the equation $H_A(\xi:\eta)=0$ determines the quadric of
all $d$-tuples $A$ whose Hessian vanishes at $(\xi:\eta)$. We'll use
freely a few facts about the Hessians, whose verification will
generally be left to the reader; the most important is the following,
which we want to highlight:

\proclaim{Lemma 3.1} The orbit of the $d$-fold point in $\P^d$ is cut out
scheme-theo\-retically by the equations $H_A(\xi:\eta)=0$,
$(\xi:\eta)\in\P^1$.\endproclaim
\demo{Proof} It is clear that the orbit of the $d$-fold point is
contained in all $H_A(\xi:\eta)$ for each $(\xi:\eta)$, so it's enough
to show that the quadrics $H_A(\xi:\eta)$ cut out the orbit at the
$d$-tuple $x^d=0$. Now the tangent space to $H_A(\xi:\eta)$ at
$(1:0:\dots:0)$ is
$$\sum_{i=0}^d i(i-1)a_i\xi^{2d-i-2}\eta^{i-2}=0\quad,$$
so the intersection of the tangent spaces at $(1:\dots:0)$ is given by
$a_2=\dots=a_d=0$, the tangent space to the orbit.\qed\enddemo

To evaluate the multiplicity of the orbit closure of a $d$-tuple at
points of its boundary, we use the techniques of \cite{Fulton},
Chapter 4: the multiplicity of a variety $Y$ along an irreducible
subvariety $X$ is the coefficient of $[X]$ in the Segre class $s(X,Y)$
of $X$ in $Y$ (\cite{Fulton}, \S 4.3), and Segre classes behave well
with respect to proper maps (\cite{Fulton}, \S 4.2). For each
component of the boundary of an orbit closure, we'll pull-back
equations for the component (essentially provided by the above lemma)
to the varieties constructed in the degree computations. Computing the
relevant term in the Segre class will be manageable on these varieties
as they are non-singular. A push-forward will then give the Segre
class in the orbit closure, and compute the multiplicity.

The boundary of the orbit closure of a $d$-tuple supported on a pair
of points consists just of the orbit of a $d$-fold point.

\proclaim{Proposition 3.2} ($s=2$) If $r\ne d/2$, the orbit closure of a
$d$-tuple consisting of one $r$-fold point and one $(d-r)$-fold point
has multiplicity 2 along its boundary. If $r=d/2$, this orbit closure
is non-singular.\endproclaim

\demo{Proof} Pull back all equations $H_A(\xi:\eta)=0$ via the map
$\P^1\times\P^1 @>>> \P^d$ considered in Proposition (1.1).
With the notations of \S1, $H_A(\xi:\eta)$ pulls back to
$$(a_1 b_0-a_0 b_1)^2(d-1)(d-r)r(a_1\xi-a_0\eta)^{2r-2}(b_1 \xi-b_0
\eta)^{2(d-r)-2}\quad;$$
as $(\xi:\eta)$ varies in $\P^1$ we see that the equations of the
orbit of the $d$-tuple pull back to the square of the equation of the
diagonal in $\P^1\times\P^1$. The diagonal maps isomorphically onto
the orbit of the $d$-fold point, and the map from $\P^1\times\P^1$ to
the orbit closure has degree 1 if $r\ne d/2$: thus, pushing forward to
$\P^d$, it follows that the first term in the Segre class of the orbit
of the $d$-fold point in the orbit closure is {\sl twice\/} the class
of the orbit. The first assertion follows. If $r=d/2$, the map from
$\P^1\times\P^1$ to the orbit closure has degree 2: thus the first
term in the Segre class is the orbit of the $d$-fold point, with
coefficient $2/2=1$. So the orbit closure is non-singular in this
case, as already observed earlier.\qed\enddemo\vskip 12pt

$s\ge 3$. If the $d$-tuple consists of at least 3 distinct points,
then its stabilizer in $\PGL(2)$ is finite, so its orbit closure is a
{\sl threefold\/} in $\P^d$. We have seen in \S 2 that the boundary of
the orbit of a $d$-tuple consists of the union of the 1-dimensional
orbit of $x^d$ and the 2-dimensional orbits of $x^ry^{d-r}$, for all
$r$ that appear as the multiplicity of a point in the $d$-tuple.

We call `premultiplicity' the product of the multiplicity of the orbit
closure of a $d$-tuple $C$ (with $s\ge 3$) by the order of its
stabilizer. Given $C$, consider its Hessian $H_C$, this time
specifically as a degree-$(2d-4)$ form on $\P^1$, and thus as a
$(2d-4)$-tuple determined by $C$. An important role is going to be
played by the points of this $(2d-4)$-tuple that lie {\sl away\/} from
$C$. We state the results first:

\proclaim{Proposition 3.3} The premultiplicity of the orbit closure of $C$
along the orbit of the $d$-fold point is
$$\sum k_i^2+4s-8\quad,$$
where the summation runs over all zeros of the Hessian $H_C$ external
to the $d$-tuple, and the $k_i$ denote the multiplicity of $H_C$ at
such points.\endproclaim

For example, suppose the Hessian is simple at all points external to
$C$; since the Hessian has degree $2d-4$, and each point with
multiplicity $r$ on $C$ contributes precisely a $(2r-2)$-fold point to
the Hessian, we find that in this case $H_C$ has exactly $2s-4$ simple
points outside of $C$, so the premultiplicity along the orbit of the
$d$-fold point must be
$$(2s-4)+(4s-8)=6(s-2)\quad.$$
In particular, the orbit closure of the general $d$-tuple, $d\ge 5$,
has multiplicity $6(d-2)$ along this orbit.

Next for the 2-dimensional components of the boundary. For every point
$p$ of $C$ of multiplicity $r$, denote by $C_p$ the residual
$(d-r)$-tuple to $p$ in $C$. In this case it matters whether the point
$p$ of $C$ is a point of the Hessian {\sl of its residual\/} $C_p$ in
$C$ (thus automatically external to $C_p$ !).

As seen in \S 2, $p$ contributes to the boundary of the orbit closure
of $C$ by the orbit of $x^ry^{d-r}$.  The next result may be seen as a
refinement of that statement:

\proclaim{Proposition 3.4} Each $r$-fold point $p$ of the $d$-tuple
contributes to the premultiplicity of the orbit closure along the
orbit of $x^ry^{d-r}$ by
$$2+\text{mult\. of $p$ in } H_{C_p}$$
if $r\ne d/2$, and
$$4+2\,(\text{mult\. of $p$ in } H_{C_p})$$
if $r=d/2$.\endproclaim

So the orbit closure of the general $d$-tuple has multiplicity $2d$
along its only boundary component (i.e., the orbit of $xy^{d-1}$),
for $d\ge5$.

\demo{Proofs} For the first computation (multiplicity along orbit of the
$d$-fold point), every point $(\xi:\eta)$ in $\P^1$ gives one equation
for the orbit of the $d$-fold point in $\P^d$, i.e\. $H_A(\xi:\eta)=0$
(see Lemma (3.1)). Now if $\phy\in\P^3$, the Hessian of the translate by
$\phy$ is given by
$$H_{A\circ\phy}=(\text{det}\phy)^2\,H_A\circ\phy\quad:$$
therefore each of the above equations for the orbit of $x^d$
pulls-back in $\P^3$ to the square of the equation of the locus $D$ of
rank-1 matrices, times the equation of the point-condition in $\P^3$
{\sl relative to the Hessian of the $d$-tuple.\/} As seen in \S 1,
point-conditions are separated above the base locus by the blow-up
resolving the rational map determined by the $d$-tuple; equations for
the inverse image of the orbit of $x^d$ in the blow-up are therefore
$$\Til D^2\,\Til H(\xi:\eta)\quad,\quad (\xi:\eta)\in\P^1$$
where $\Til D$ is the equation for the proper transform of $D$, and
$\Til H(\xi:\eta)$ is the point-condition in the blow-up {\sl relative
to the points in the Hessian not contained in the $d$-tuple\/.} The
scheme-theoretic inverse image consists then of a non-reduced
scheme supported on the proper transform of the locus of rank-1
matrices, with length 2 over the support, and embedded components
along pencils of matrices whose image is a point of the Hessian {\sl
not\/} contained in the $d$-tuple. Each of the embedded pencils maps
isomorphically to the 1-dimensional orbit of $x^d$, and contributes to
the 1-dimensional term of the Segre class of the scheme by the square
of the multiplicity of the corresponding point in the Hessian; this
accounts for the $\sum k_i^2$ term in the formula.  It remains
therefore to be seen that the proper transform $\Til D$ of the locus
of rank-1 matrices accounts for the term $4s-8$ in the
premultiplicity. Now we claim that all we have to check is that $\Til
D^2$ pushes forward to $(2-s)$ times the class of the orbit of $x^d$:
indeed, it will follow that the contribution of $\Til D$ to the
1-dimensional term of the Segre class (i.e., $-(2\Til D)^2$) pushes
forward in $\P^d$ to $(4s-8)$ times the class of the orbit of $x^d$,
and we will be done. Now a straightforward computation shows that the
push-forward of $\Til D^2$ is the push-forward {\sl from\/} $\P^3$ of
$D^2$ {\sl minus\/} the $s$ lines of the base locus (which map
isomorphically to the orbit of $x^d$). Finally, $D^2$ consists, as a
class on the quadric $D$, of 2 lines of each ruling, and the ruling
parametrizing matrices with given kernel pushes forward to 0 in
$\P^d$; so the push-forward is indeed $2-s$ times the orbit, as
needed.\vskip 12pt

For the second statement (the multiplicity along the orbit of
$x^ry^{d-r}$), suppose $p$ is a point of multiplicity $r$ in the
$d$-tuple, and factor the map $\P^3 \dashrightarrow \P^d$ through
$$\P^3 \dashrightarrow \P^1\times\P^{d-r} @>>>\P^d\quad,$$ where
$\P^3$ maps to each factor $\P^1$ and $\P^{d-r}$ as usual, by
extending the action of $\PGL(2)$ on the $r$-fold point $p$ and its
residual $(d-r)$-tuple $C_p$ respectively; the orbit closure of this
point $(p,C_p)$ in $\P^1\times\P^{d-r}$ maps surjectively to
the orbit closure of the $d$-tuple in $\P^d$. The point is that the
map $\P^1\times\P^{d-r} @>>>\P^d$ is an immersion at every point
$(p,(d-r)q)$ {\sl if\/} $p\ne q$; moreover, in this case the inverse
image of $rp+(d-r)q$ consists of precisely $(p,(d-r)q)$ if $r\ne d/2$,
and of the two points $(p,(d-r)q)$ and $(q,rp)$ if $r=d/2$. Thus we
only have to show that the premultiplicity of the orbit closure of
$(p,C_p)$ in $\P^1\times\P^{d-r}$ is $2+$mult\. of $p$ in the
Hessian of $C_p$.

For this, we observe that equations for the set of points in
$\P^1\times\P^{d-r}$ of type $(p,(d-r)q)$ are (again by Lemma (3.1))
given by $H_A(\xi:\eta)=0$, where now the Hessian is taken for
$A\in\P^{d-r}$. Pulling back to $\P^3$, and recalling again that the
Hessian of a translate is the translate of the Hessian multiplied by
the square of the determinant of the translation, we find that
equations in $\P^3$ for the inverse image of the locus of pairs
$(p,(d-r)q)$ are
$$(det\phy)^2\,H_{C_p}(\phy(\xi:\eta))=0.$$

Now blow-up $\P^3$ as usual, and study it over the pencil of all
$\phy$ whose image is the $r$-fold point $p$ of the $d$-tuple. By
arguing as in \S1, one sees that the blow-up resolves the map $\P^3
\dashrightarrow \P^1\times\P^{d-r}$; pulling back the above equation
to the blow-up, we find that (near the pencil) the inverse image of
the locus of pairs $(p,(d-r)q)$ is supported on the proper transform
of the determinant hypersurface (with length 2), and on the component
of the exceptional divisor over the pencil (with length $2+$mult\. of
$p$ in $H_{C_p}$). Now pairs $(p,(d-r)q)$ {\sl with $p\ne q$\/} don't
come from the determinant hypersurface (which maps to $d$-fold points
only), so the premultiplicity equals the length of the part supported
on the exceptional divisor, and this concludes the proof of the last
claim.\qed\enddemo

\heading \S 4. Smooth orbit closures and more. \endheading

The results of \S3, together with a description of the finite subgroups
of \PGL(2) (see [{\bf Weber}], \S\S67-77), allow us to give an immediate
classification of the smooth \PGL(2)-orbit closures.

First we present the following lemma, some instances of which appeared
already above. Its proof may be left to the reader.

\proclaim{Lemma 4.1} The map $\P^d \to\P^{md}$, $f\mapsto f^m$ is an embedding.
\endproclaim

If the $d$@-tuple corresponding to $f$ is supported on $s\ge3$ points,
the orbit closure of $f^m$ has degree equal to $m^3$ times the degree
of the orbit closure of $f$ (for example by Proposition 1.3), whereas
the multiplicities along corresponding boundary components are equal.

Because of the lemma, in the remainder of this section we will only
consider $d$@-tuples for which the g.c.d.~of the multiplicities of the
$s$ points equals one. We will also assume that $s\ge3$; recall that
the orbit (closure) of $x^d$ is smooth and that the orbit closure of
$x^ry^{d-r}$ is smooth if and only if $d=2r$.

With these assumptions, we have:

\proclaim{Proposition 4.2} The smooth 3-dimensional \PGL(2)-orbit
closures are:
\roster
\item the orbit closure of $x^3+y^3$, with stabilizer $D_3=S_3$;
\item the orbit closure of $x^4+xy^3$, with stabilizer $A_4$;
\item the orbit closure of $x^5y-xy^5$, with stabilizer $S_4$;
\item the orbit closure of $x^{11}y+11x^6y^6-xy^{11}$, with stabilizer
$A_5$.
\endroster
\endproclaim

\demo{Proof} The orbit closure of a $d$@-tuple $f$ is smooth if and only if its
multiplicity along the orbit of $x^d$ equals one, i.e., the premultiplicity
along that orbit equals the order of the stabilizer of $f$. {From}
Proposition (3.3),
this premultiplicity equals $\sum k_i^2+4s-8$, where the $k_i$ are the
multiplicities of the points of the Hessian of $f$ external to $f$. Counted
with multiplicity, there are $2s-4$ such points (i.e., $\sum k_i = 2s-4$), so
the premultiplicity is $\ge6(s-2)$.

Assuming that $f$ has smooth orbit closure, it follows that the order of its
stabilizer is $\ge6(s-2)$. In particular, its stabilizer is non-trivial.
It now suffices to consider the action of the finite subgroups $G$ of \PGL(2)
on $\P^1$ and the orbits of points with non-trivial stabilizer. Following
[{\bf Weber}], \S68, we list these groups and the lengths of the special
orbits:
\roster
\item[0] $G=C_n$; lengths 1, 1;
\item $G=D_n$; lengths 2, $n$, $n$;
\item $G=A_4$; lengths 4, 4, 6;
\item $G=S_4$; lengths 6, 8, 12;
\item $G=A_5$; lengths 12, 20, 30.
\endroster
Determining the $d$@-tuples $f$ with smooth orbit closure is now an easy
matter:
\roster
\item[0] Assume $\text{Stab}(f)=C_n$. Then $n\ge6(s-2)>s$. It follows that
$f$ is supported on one or two points, a contradiction.
\item Assume $\text{Stab}(f)=D_n$. Then $2n\ge6(s-2)$ so $n\ge3(s-2)\ge s$.
Again, if $n>s$ it follows that $s=2$, a contradiction; so we get $n=s=3$ and
$\text{Stab}(f)=D_3=S_3$. Clearly the multiplicities of the 3 points are
all equal, thus by our assumption they are all one. So this is the
orbit closure of $x^3+y^3$, which is $\P^3$. Of course smoothness also follows
from considering the Hessian of $f$.
\item Assume $\text{Stab}(f)=A_4$. Then $12\ge6(s-2)$ so $s\le4$. It follows
that $s=4$ and that all multiplicities are equal (to one). This is the orbit
closure of $x^4+xy^3$; computing the Hessian, we see that it is indeed
smooth.
\item Assume $\text{Stab}(f)=S_4$. Then $24\ge6(s-2)$ so $s\le6$. It follows
that $s=6$ and that all multiplicities are equal to one. This is the orbit
closure of $x^5y-xy^5$, which is indeed smooth, as its Hessian has simple
zeros.
\item Assume $\text{Stab}(f)=A_5$. Then $60\ge6(s-2)$ so $s\le12$. It follows
that $s=12$ and that all multiplicities are equal to one. This is the orbit
closure of $x^{11}y+11x^6y^6-xy^{11}$ ([{\bf Weber}], \S74). It is smooth
as its Hessian has 20 simple zeros.\qed
\endroster
\enddemo

It turns out that it is also possible to classify the orbit closures
that are smooth in codimension one. The answer is particularly pretty
in case the multiplicities of the $s$ points of the $d$@-tuple are all
equal. In that case we may and will assume that they are all equal to
one, so that $d=s$; call such a $d$@-tuple {\sl simple\/}. Note that
the orbit closure of a simple $d$@-tuple has at most one boundary
component.

\proclaim{Proposition 4.3} The orbit closure of a simple $d$@-tuple $f$ is
smooth in codimension one if and only if $f$ is a special orbit for
the action of a finite subgroup $G$ of \PGL(2) on $\P^1$ (i.e., $f$ is
an orbit of length smaller than the order of $G$).\endproclaim

\demo{Proof} Let $f$ be a simple $d$@-tuple (so $d=s$). If $d=1$
(resp.~2) the orbit closure of $f$ is smooth; take $G=C_n$ (resp.
$D_n$) for an $n\ge2$.  So we assume $d\ge3$. {From} Proposition (3.4),
the premultiplicity of the orbit closure of $f$ along its only
boundary component equals $\sum (2+\text{mult\. of $p$ in } H_{C_p})$,
where the summation runs over the $d$ points $p$ of $f$. Assuming that
the orbit closure of $f$ is smooth in codimension one, it follows that
the stabilizer of $f$ has order $\ge 2d$. The ``only if" part of the
proposition follows. It remains to check that the orbit closures of
the special orbits are indeed smooth in codimension one.  This is an
easy verification (see below).\qed\enddemo

It is perhaps worthwhile to remark that the proposition above seems to
constitute an answer to the question raised in [{\bf M-U}], Remark
(3.6): the \PGL(2)@-orbit closures of special $G$@-orbits
($G\subset\PGL(2)$ finite) may be characterized as the orbit closures
of simple $d$@-tuples that are smooth in codimension one.

The general case is somewhat harder. Let $f$ be a $d$@-tuple supported
on $s\ge3$ points, and assume that the orbit closure of $f$ is smooth
in codimension one. Suppose that there are $s_a$ points with
multiplicity $a$.  Then the stabilizer of $f$ has order at least
$2s_a$. We conclude that $f$ is supported on the special orbits for
the action of its stabilizer $G$ on $\P^1$. Clearly $G$ is not cyclic,
so there are 3 such orbits. Call them $A$, $B$ and $C$, and write
$f=A^aB^bC^c$ with $a$, $b$ and $c$ positive integers. Call
$A$@-multiplicity the contribution of the points of $A$ to the
multiplicity of the orbit closure of $f$ along the orbit of
$x^ay^{d-a}$.  By Proposition (3.4), this equals
$$\frac{d_A(2+\text{mult. of $p$ in the Hessian of }A_p^aB^bC^c)}
{\text{order of }G}
$$
where $d_A$ is the degree of $A$, $p$ a point of $A$ and $A_p$ the residual
$(d_A-1)$@-tuple.

Similarly we define the $B$@-multiplicity and the $C$@-multiplicity. The
following result is an immediate consequence.

\proclaim{Proposition 4.4} Let $G$ be a finite, non-cyclic subgroup of
\PGL(2). Denote by $A$, $B$ and $C$ the three special orbits for the
action of $G$ on $\P^1$. Let $f=A^aB^bC^c$, with $a$, $b$ and $c$
positive integers. Assume that $G$ is the \PGL(2)@-stabilizer of $f$.
The \PGL(2)@-orbit closure of $f$ is smooth in codimension one if and
only if $a$, $b$ and $c$ are mutually distinct and the
$A$@-multiplicity, the $B$@-multiplicity and the $C$@-multiplicity are
equal to one.\endproclaim

When one or two of $a$, $b$ and $c$ are zero, the proposition remains true,
{\sl mutatis mutandis\/}.

Computing the multiplicity of the Hessian at $p$ becomes simpler when
one chooses the right coordinates. Namely, $p$ is one of
the two fixed points of an element of $G$ (of order $m=(\text{order of }
G)/d_A$). Choose coordinates $x$, $y$ so that $p$ and the other fixed
point are given by $x=0$ and $y=0$ respectively. Writing out
$A_p$, $B$ and $C$ in these coordinates, we see that only powers of
$x^m$ occur:
$$\align
A_p &=y^{d_A-1}+A_1y^{d_A-1-m}x^m+A_2y^{d_A-1-2m}x^{2m}+\dots,\\
B &=y^{d_B}+B_1y^{d_B-m}x^m+B_2y^{d_B-2m}x^{2m}+\dots,\\
C &=y^{d_C}+C_1y^{d_C-m}x^m+C_2y^{d_C-2m}x^{2m}+\dots\,.
\endalign $$
Now one immediately checks that the multiplicity of the Hessian
of $A_p^aB^bC^c$ at $p$ is $m-2$ when
$$A_1a+B_1b+C_1c\neq0;$$
that it is $2m-2$ when
$$\align
A_1a+B_1b+C_1c &=0 \text{ and}\\
(A_1^2-2A_2)a+(B_1^2-2B_2)b+(C_1^2-2C_2)c &\neq0,
\endalign$$
etc. Thus the $A$@-multiplicity is 1, 2, \dots, correspondingly.

Finally we list for each of the finite, non-cyclic subgroups
$G$ of \PGL(2) the special orbits and the relevant equations:
\roster
\item $G=D_n$: $A=xy$, $B=x^n+y^n$, $C=x^n-y^n$; the
$A$@-multiplicity is 1 iff $$b\neq c;$$ the $B$@-multiplicity is
1 iff $$-a+\frac{(n-1)(n-2)}{6}b+\frac{n(n-1)}{2}c\neq0;$$ the
$C$@-multiplicity is 1 iff $$-a+\frac{n(n-1)}{2}b+
\frac{(n-1)(n-2)}{6}c\neq0.$$
\item $G=A_4$: $A=x^4+2\sqrt{-3}x^2y^2+y^4$, $B=x^4-2\sqrt{-3}x^2y^2
+y^4$, $C=x^5y-xy^5$; the $A$@-multiplicity is 1 iff $$a-8b+20c\neq0,
$$ otherwise it is 2;
the $B$@-multiplicity is 1 iff $$8a-b-20c\neq0,$$ otherwise it is 2;
the $C$@-multiplicity
is 1 if $$a\neq b;$$ it is 2 when $a=b$ (unless $c=14a=14b$, in
which case it is 4); note however that when $a=b$ the stabilizer
is $S_4$, so the actual multiplicities are 1, respectively 2 (see
also below).
\item $G=S_4$: $A=x^5y-xy^5$, $B=x^8+14x^4y^4+y^8$, $C=
x^{12}-33x^8y^4-33x^4y^8+y^{12}$; the $A$@-multiplicity is 1 iff
$$a-14b+33c\neq0,$$ otherwise it is 2; the $B$@-multiplicity is 1 iff
$$20a-7b-88c\neq0,$$ otherwise it is 2; the $C$@-multiplicity is 1
iff $$45a-84b-11c\neq0,$$ it is 2 when $45a-84b-11c=0$, unless
$(a,b,c)\sim(5852,561,19656)$, in which case it is 3.
\item $G=A_5$: $A=x^{11}y+11x^6y^6-xy^{11}$, $$\align
B &=x^{20}-228x^{15}y^5+494x^{10}y^{10}+228x^5y^{15}+y^{20},\\
C
&=x^{30}+522x^{25}y^5-10005x^{20}y^{10}-10005x^{10}y^{20}-522x^5y^{25}+y^{30};
\endalign $$
the $A$@-multiplicity is 1 iff $$11a-228b+522c\neq0,$$ otherwise it
is 2; the $B$@-multiplicity is 1 iff $$88a-57b-580c\neq0,$$
otherwise it is 2; the $C$@-multiplicity is 1 iff
$$99a-285b-58c\neq0,$$ it is 2 otherwise, unless
$(a,b,c)\sim(26864005,431607,43733250)$, in which case it is 3.

\endroster

\enddocument